\documentclass[a4paper,11pt]{amsart}
\usepackage{fancyhdr}
\usepackage{amsmath}
\usepackage{dsfont}
\usepackage{hyperref}
\usepackage[mathscr]{eucal}
\usepackage[cp1251]{inputenc}
\usepackage[english]{babel}
\usepackage{enumerate,float,indentfirst}
\usepackage{graphicx}
\usepackage{xcolor}
\usepackage{latexsym,a4,mathrsfs,amsthm,amsmath,amssymb,url}
\usepackage{amsfonts}
\usepackage{amssymb}
\usepackage{caption}

\usepackage{tikz}
\usepackage{pgfplots}
\usetikzlibrary{plotmarks}

\usetikzlibrary{decorations.pathreplacing}

\usepackage{stackengine}

\numberwithin{equation}{section}
\newtheorem{hypothesis}{Hypothesis}
\newtheorem{problem}{Problem}

\begin{document}

\vspace{0in}

\title[Submatrices with the best-bounded inverses]{\bf Submatrices with the best-bounded inverses: Studying $\mathds{R}^{n \times 2}$ and $\mathds{C}^{n \times 2}$}

\author[Yu. Nesterenko]{Yuri Nesterenko}
\email{y\_nesterenko@mail.ru}

\begin{abstract}
In both real and complex cases, we establish the connection of the problem about $2$-dimensional linear subspaces the most deviating from the coordinate ones with one simply formulated optimization problem for isoperimetric polygons in Euclidean spaces.

This study thereby provides a new geometrical point of view on the $2$-dimensional case of the problem formulated by Goreinov, Tyrtyshnikov and Zamarashkin \cite{GTZ1997}, and at the same time presents a new application of the results by Hausmann and Knutson \cite{HK1997}.
\end{abstract}

\maketitle

\thispagestyle{empty}
\vspace{-5truemm}
\section{Introduction}
In this paper we study the $2$-dimensional case of the problem about linear subspaces the most deviating from the coordinate ones.

\begin{problem}\label{p1}
Given $n > k > 0$ and $\mathds{F} = \mathds{R} \, (\mathds{C})$, what is the largest possible deviation of a $k$-dimensional subspace of $\mathds{F}^n$ from all the $k$-dimensional coordinate subspaces, if the deviation is measured as the largest principal angle?
\end{problem}

Initially formulated by Goreinov, Tyrtyshnikov and Zamarashkin \cite{GTZ1997} as the hypothesis about matrices and their submatrices, it concerned the real case.

\begin{hypothesis}\label{f1}
For every $n > k > 0$ and an arbitrary real $n \times k$ matrix with orthonormal columns, such a $k \times k$ submatrix exists that the spectral norm of its inverse does not exceed $\sqrt{n}$.
\end{hypothesis}

According to this hypothesis, the answer for the real version of the Problem \ref{p1} is bounded above by $\arccos(1 / \sqrt{n})$. Moreover, various numerical experiments allow us to suggest that this estimate is sharp, and the corresponding equality is achieved for a finite number of subspaces.

At the same time, the problem remains open for all $n$ and $k$ such that $1 < k < n - 1$, except $n = 4, \, k = 2, \, \mathds{F} = \mathds{R}$ (see \cite{Nesterenko2023}).

In the next section we study the case $k = 2, \, \mathds{F} = \mathds{R}$ in more detail by establishing its connection with one simply formulated optimization problem for isoperimetric polygons on the Euclidean plane.

Further, we do the same for the complex case\footnote{For our knowledge, there is no a ready hypothesis for it. The estimate from the real case is violated. The counterexample is available at the end of Section \ref{secc}. } by establishing the similar connection with isoperimetric polygons in the $3$-dimensional Euclidean space.

\section{$\mathds{R}^{n \times 2}$}

In this section we consider matrix $A \in \mathds{R}^{n \times 2}$ such that $A^T A = I$. Let's reformulate this property in terms of $A$'s rows.

Since both singular values of $A$ equal $1$, the sum of squared projections of an arbitrary unit vector $x \in \mathds{R}^2$ onto the rows of $A$ equals $1$. Therefore,
\begin{equation*}
\sum_{i = 1}^{n} r_i^2 \cos^2{(\varphi - \varphi_i)} = 1, \quad \forall \varphi \in [0, 2\pi),
\end{equation*}
where $r_i$ and $\varphi_i$ are polar coordinates of $A$'s rows, and $\varphi$ is the polar angle of vector $x$.

Since $\sum_{i = 1}^{n} r_i^2 = 2$, the previous equation yields
\begin{equation*}
\sum_{i = 1}^{n} r_i^2 \cos{(2\varphi - 2\varphi_i)} = 0, \quad \forall \varphi \in [0, 2\pi),
\end{equation*}
or 
\begin{equation*}
\sum_{i = 1}^{n} ( a_i, y ) = 0, \quad \forall y \in \mathds{R}^2: \, |y| = 1,
\end{equation*}
where $a_1, \ldots, a_n \in \mathds{R}^2$ are complex-like squares of the rows of $A$, and $y$ is the similar square of $x$.

Therefore, vectors $a_1, \ldots, a_n$ constitute a plane polygon of perimeter $2$
\begin{equation*}
\sum_{i = 1}^{n} a_i = \theta, \quad \sum_{i = 1}^{n} | a_i | = 2.
\end{equation*}

Note that if we transform matrix $A$ multiplying it by a $2 \times 2$ orthogonal matrix, the corresponding polygon will undergo an orthogonal transformation too.

So far, we have got only an alternative derivation of the fact given in \cite{HK1997}.

Let's apply the same technique to the estimate from the hypothesis we study. It states that there exist such $i, j \in {1, \ldots, n}$ that
\begin{equation*}
\min_{\varphi} \, ( r_i^2 \cos^2{(\varphi - \varphi_i)} + r_j^2 \cos^2{(\varphi - \varphi_j)} ) \geq \frac{1}{n}.
\end{equation*}

After reduction we obtain
\begin{equation*}
\min_{\varphi} \, ( r_i^2 \cos{(2\varphi - 2\varphi_i)} + r_j^2 \cos{(2\varphi - 2\varphi_j)} ) \geq \frac{2}{n} - r_i^2 - r_j^2.
\end{equation*}

In terms of vectors $a_1, \ldots, a_n$ and $y$, it means that
\begin{equation*}
\min_{|y| = 1} ( a_i + a_j, y ) \geq \frac{2}{n} - |a_i| - |a_j|,
\end{equation*}
which is equivalent to the inequality
\begin{equation*}
|a_i| + |a_j| - |a_i + a_j| \geq \frac{2}{n}.
\end{equation*}

We are ready now to put forward a self-contained hypothesis which is equivalent to original one with $k = 2$, while also providing it with the equality criteria. Note that vectors $a_1, \ldots, a_n$ below are renormalized.

\begin{hypothesis}\label{hr}
For arbitrary vectors $a_1, \ldots, a_n, \, n \geq 3$ forming a unit-perimeter planar polygon
\begin{equation*}
\begin{split}
1) \, & \max_{i,j} (|a_i| + |a_j| - |a_i + a_j|) \geq \frac{1}{n}, \\
2) \, & \text{equality holds} \iff \{ a_1, \ldots, a_n \} = \{ a, b, c \}, \, \text{where} \\
& |a| = \frac{1}{4n} + \frac{1}{4p}, \; p = \# \{ i = 1, \ldots, n: \, a_i = a \} \geq 1, \\
& |b| = \frac{1}{4n} + \frac{1}{4q}, \; q = \# \{ i = 1, \ldots, n: \, a_i = b \} \geq 1, \\
& |c| = \frac{1}{4n} + \frac{1}{4r}, \; r = \# \{ i = 1, \ldots, n: \, a_i = c \} \geq 1.
\end{split}
\end{equation*}
\end{hypothesis} 

The following picture illustrates the point 2. The order of vectors $a_1, \ldots, a_n$ is chosen in such a way that the polygon they form is convex.

\begin{center}
\begin{figure}[h]
\begin{tikzpicture}

\def \d {0.04cm}
\def \r {0.03cm}

\begin{scope}[shift={(160 * \d,-160 * \d)}]

\draw[fill=black][black] (-20 * \d, 80 * \d) circle (\r);

\draw [decorate,decoration={brace,amplitude=5pt,mirror,raise=4ex}]
  (-20 * \d, 80 * \d) -- (-50 * \d, -10 * \d) node[midway]{};

\draw[fill=white][white] (-56 * \d, 42 * \d) node[anchor=east,text=black]{$p$};

\draw[fill=white][white] (-25 * \d, 67 * \d) node[anchor=east,text=black]{$a$};

\draw[-stealth][black] (-20 * \d, 80 * \d) -- (-30 * \d, 50 * \d);

\draw[fill=black][black] (-30 * \d, 50 * \d) circle (\r);

\draw[black] (-28 * \d, 56 * \d) -- (-30 * \d, 50 * \d);

\draw[fill=black][black] (-33 * \d, 41 * \d) circle (0.25 * \r);
\draw[fill=black][black] (-35 * \d, 35 * \d) circle (0.25 * \r);
\draw[fill=black][black] (-37 * \d, 29 * \d) circle (0.25 * \r);

\draw[fill=black][black] (-40 * \d, 20 * \d) circle (\r);

\draw[-stealth][black] (-40 * \d, 20 * \d) --  (-50 * \d, -10 * \d);

\draw[fill=white][white] (-45 * \d, 7 * \d) node[anchor=east,text=black]{$a$};

\draw[fill=black][black] (-50 * \d, -10 * \d) circle (\r);

\draw [decorate,decoration={brace,amplitude=5pt,mirror,raise=4ex}]
  (-50 * \d, -10 * \d) -- (100 * \d, -10 * \d) node[midway]{};

\draw[fill=white][white] (25 * \d, -32 * \d) node[anchor=north,text=black]{$q$};

\draw[-stealth][black] (-50 * \d, -10 * \d) --  (-10 * \d, -10 * \d);

\draw[fill=white][white] (-30 * \d, -12 * \d) node[anchor=north,text=black]{$b$};

\draw[fill=black][black] (-10 * \d, -10 * \d) circle (\r);


\draw[fill=black][black] (19 * \d, -10 * \d) circle (0.25 * \r);
\draw[fill=black][black] (25 * \d, -10 * \d) circle (0.25 * \r);
\draw[fill=black][black] (31 * \d, -10 * \d) circle (0.25 * \r);


\draw[fill=black][black] (60 * \d, -10 * \d) circle (\r);

\draw[-stealth][black] (60 * \d, -10 * \d) --  (100 * \d, -10 * \d);

\draw[fill=white][white] (80 * \d, -12 * \d) node[anchor=north,text=black]{$b$};

\draw[fill=black][black] (100 * \d, -10 * \d) circle (\r);

\draw [decorate,decoration={brace,amplitude=5pt,mirror,raise=4ex}]
  (100 * \d, -10 * \d) -- (-20 * \d, 80 * \d) node[midway]{};

\draw[fill=white][white] (52 * \d, 51 * \d) node[anchor=south west,text=black]{$r$};

\draw[-stealth][black] (100 * \d, -10 * \d) --  (55 * \d, 23.75 * \d);

\draw[fill=white][white] (77.5 * \d, 6.875 * \d) node[anchor=south west,text=black]{$c$};

\draw[fill=black][black] (55 * \d, 23.75 * \d) circle (\r);

\draw[fill=black][black] (34 * \d, 39 * \d) circle (0.25 * \r);
\draw[fill=black][black] (40 * \d, 35 * \d) circle (0.25 * \r);
\draw[fill=black][black] (46 * \d, 31 * \d) circle (0.25 * \r);

\draw[fill=black][black] (25 * \d, 46.25 * \d) circle (\r);

\draw[-stealth][black] (25 * \d, 46.25 * \d) --  (-20 * \d, 80 * \d);

\draw[fill=white][white] (2.5 * \d, 63.125 * \d) node[anchor=south west,text=black]{$c$};

\end{scope}

\end{tikzpicture}
\caption{ Polygons described in Hypothesis \ref{hr}. }
\label{triangle}
\end{figure}
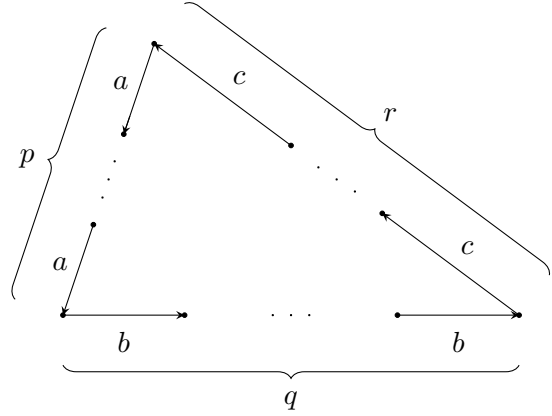
\end{center}

It can be shown using the cosine rule that regardless the partition of vectors $a_1, \ldots, a_n$ into three groups of equal to $a$, $b$ and $c$
\begin{equation*}
|a| + |b| - |a + b| = |a| + |c| - |a + c|
= |b| + |c| - |b + c| = \frac{1}{n}.
\end{equation*}

Even though the whole statement is still waiting to be proven\footnote{Or disproven.} for $n > 4$, the presented reformulation provides us with a new intuition about the original hypothesis.

\section{$\mathds{C}^{n \times 2}$}\label{secc}

In this section $A$ will be a complex $n \times 2$ matrix such that $A^H A = I$. Let's denote its rows as $(u_i, v_i) \in \mathds{C}^2, \, i = 1, \ldots, n$.

Since both singular values of $A$ equal $1$, taking an arbitrary unit vector $(u, v) \in \mathds{C}^2$ and computing the sum of its squared projections onto the rows of $A$ yields~$1$. Therefore,
\begin{equation}\label{expr_c0}
\sum_{i = 1}^{n} | u_i \bar{u} + v_i \bar{v} |^2 = 1, \quad \forall u, v \in \mathds{C}: |u|^2 + |v|^2 = 1.
\end{equation}
Since $\sum_{i = 1}^{n} |u_i|^2 = \sum_{i = 1}^{n} |v_i|^2 = 1$ and $|u|^2 + |v|^2 = 1$, we can rewrite (\ref{expr_c0}) as follows
\begin{equation}\label{expr_c05}
\sum_{i = 1}^{n} ( \bar{u_i} v_i u \bar{v} + u_i \bar{v_i} \bar{u} v + |u_i|^2 |u|^2 + |v_i|^2 |v|^2 - |u_i|^2 |v|^2 - |v_i|^2 |u|^2) = 0.
\end{equation}
Taking
\begin{equation*}
\begin{split}
a_i = ( \bar{u_i} v_i &+ u_i \bar{v_i}, \, i(\bar{u_i} v_i - u_i \bar{v_i}), \, |u_i|^2 - |v_i|^2 ) \in \mathds{R}^3, \; i = 1, \ldots, n,\\
&a = ( \bar{u} v + u \bar{v}, \, i(\bar{u} v - u \bar{v}), \, |u|^2 - |v|^2 ) \in \mathds{R}^3,
\end{split}
\end{equation*}
equation (\ref{expr_c05}) becomes as follows
\begin{equation*}
\sum_{i = 1}^{n} ( a_i, a ) = 0,
\end{equation*}
where $(\cdot,\cdot)$ is the standard dot product in $\mathds{R}^3$.

Since the map $(u, v) \mapsto a$ defines the well-known Hopf bundle of the unit sphere in $\mathds{C}^2$ over the unit sphere in $\mathds{R}^3$ (see \cite{Hopf1931}), we obtain $\sum_{i = 1}^{n} a_i = \theta$.

In addition, $|a_i| = |u_i|^2 + |v_i|^2, \, i = 1, \ldots, n$, and hence $\sum_{i = 1}^{n} | a_i | = 2$.

Thus, we deal again with isoperimetric polygons, but this time in $\mathds{R}^3$.

And again, if we apply a unitary transformation to the columns of $A$, the corresponding spatial polygon will undergo an orthogonal transformation (see section 3 in the paper \cite{HK1997}). 

Now, let's assume that $\arccos(B_n)$ is the answer\footnote{Further we discuss the question about its concrete value.} for the complex version of the Problem \ref{p1}. Then
\begin{equation*}
\max_{i,j} \min_{|u|^2 + |v|^2 = 1} \, \left( | u_i \bar{u} + v_i \bar{v} |^2 + | u_j \bar{u} + v_j \bar{v} |^2 \right) \geq B_n^2.
\end{equation*}

Expanding the squared modules above yields the expression
\begin{equation}\label{expr_c1}
(|u_i|^2 + |u_j|^2) |u|^2 + (|v_i|^2 + |v_j|^2) |v|^2 + 2 \operatorname{Re} (u_i \bar{v_i} + u_j \bar{v_j}) \bar{u} v.
\end{equation}

If we parameterize the unite $3$-sphere $|u|^2 + |v|^2 = 1$ as follows
\begin{equation*}
u = \cos\varphi \, e^{i \xi_1}, \, v = \sin\varphi \, e^{i \xi_2}, \, \varphi \in [0, \frac{\pi}{2}], \, \xi_1, \xi_2 \in [0, 2\pi),
\end{equation*}
and apply this substitution to the expression (\ref{expr_c1}), we obtain
\begin{equation}\label{expr_c2}
\begin{split}
\frac{1 + \cos 2\varphi}{2} (|u_i|^2 + |u_j|^2) &+ \frac{1 - \cos 2\varphi}{2} (|v_i|^2 + |v_j|^2) +\\
+ \sin 2\varphi \operatorname{Re} (u_i \bar{v_i} &+ u_j \bar{v_j}) e^{i (\xi_1 - \xi_2)}.
\end{split}
\end{equation}

Since $\sin 2\varphi \geq 0$, the minimization of (\ref{expr_c2}) over $\varphi, \xi_1, \xi_2$ can be changed to the iterated one of the form $\min_{\varphi} \min_{\xi_1, \xi_2} (\ldots)$. The inner minimum gives us the following result
\begin{equation}\label{expr_c3}
\begin{split}
&\qquad \qquad \qquad \qquad \frac{1}{2} (|u_i|^2 + |v_i|^2 + |u_j|^2 + |v_j|^2) +\\
&+ \frac{1}{2} \min_{\varphi} \left( \cos 2\varphi (|u_i|^2 - |v_i|^2 + |u_j|^2 - |v_j|^2) - \sin 2\varphi | 2 u_i \bar{v_i} + 2 u_j \bar{v_j}| \right).
\end{split}
\end{equation}

Vectors $(\cos 2\varphi, -\sin 2\varphi) \in \mathds{R}^2$ and $(|u_i|^2 - |v_i|^2 + |u_j|^2 - |v_j|^2, | 2 u_i \bar{v_i} + 2 u_j \bar{v_j}|) \in \mathds{R}^2$ lie in the upper and the lower half-plane, respectively. Therefore, the minimum of their dot product is achieved and equal to the length of the second vector with a minus sign.

Thus, the expression (\ref{expr_c3}) is equal to
\begin{equation*}
\begin{split}
&\qquad \qquad \quad \frac{1}{2} (|u_i|^2 + |v_i|^2 + |u_j|^2 + |v_j|^2) -\\
&- \frac{1}{2} \sqrt{ (|u_i|^2 - |v_i|^2 + |u_j|^2 - |v_j|^2)^2 + | 2 u_i \bar{v_i} + 2 u_j \bar{v_j}|^2 },
\end{split}
\end{equation*}
or simply
\begin{equation*}
\frac{1}{2} ( |a_i| + |a_j| - |a_i + a_j| )
\end{equation*}
in terms of the Hopf map.

Let's summarize the above. The problem about $2$-dimensional complex subspaces the most deviating from the coordinate ones with $\arccos(B_n)$ as the sharp upper bound is equivalent to the following 
\begin{problem}\label{pc}
For arbitrary vectors $a_1, \ldots, a_n, \, n \geq 3$ forming a unit-perimeter polygon\footnote{As in the previous section, vectors $a_1, \ldots, a_n$ are renormalized.} in the $3$-dimensional Euclidean space, the sharp lower bound for the expression $\max_{i,j} (|a_i| + |a_j| - |a_i + a_j|)$
is equal to $B_n^2$.
\end{problem}

Our experiments let us suggest that for $n \geq 4$ this lower bound is achieved when
\begin{equation*}
\begin{split}
& \{ a_1, \ldots, a_n \} = \{ a, b, c, d \} \; \text{and} \\
\left\lfloor \frac{n}{4} \right\rfloor \leq \, \# \{ i &= 1, \ldots, n: \, a_i = x \} \leq \left\lceil \frac{n}{4} \right\rceil, \, \forall x \in \{ a, b, c, d \}.
\end{split}
\end{equation*}

That is, in contrast to the real case, the partitions of $n$ matter, and the lower bound is achieved on the most balanced ones.

We have omitted the case $n = 3$ here, because it is equivalent to the real one and also breaks the proposed pattern.

As well, we observe that
\begin{equation*}
n B_n^2 \geq 2 - 2\sqrt{\frac{1}{3}},
\end{equation*}
where the equality holds if and only if $n \, \vdots \, 4$. The following plot illustrates the suggested behavior of $n B_n^2$.

\quad

\begin{center}
\begin{tikzpicture}
\begin{axis}
[
        height=0.45\textwidth,
        width=0.7\textwidth,
        xlabel=n,
        xtick={4,8,...,40},
        ylabel=$n B_n^2$,
        yticklabel style={%
        /pgf/number format/.cd,
        precision=3
        }
]

\addplot[mark=*,mark options={fill=white},black] coordinates {
(4, 0.845299462)
(5, 0.85509    )
(6, 0.857142   )
(7, 0.853503   )
(8, 0.845299462)
(9, 0.8486118  )
(10, 0.849473   )
(11, 0.848276   )
(12, 0.845299462)
(13, 0.8469526  )
(14, 0.8474172  )
(15, 0.846831   )
(16, 0.845299462)
(17, 0.8462872  )
(18, 0.8465778  )
(19, 0.8462315  )
(20, 0.845299462)
(21, 0.8459556  )
(22, 0.846153   )
(23, 0.8459262  )
(24, 0.845299462)
(25, 0.8457675  )
(26, 0.84591    )
(27, 0.845748   )
(28, 0.845299462)
(29, 0.8456487  )
(30, 0.845757   )
(31, 0.8456366  )
(32, 0.845299462)
(33, 0.8455722  )
(34, 0.8456582  )
(35, 0.845565   )
(36, 0.845299462)
(37, 0.8455166  )
(38, 0.8455836  )
(39, 0.8455122  )
(40, 0.845299462)};

\addplot[black, dotted] coordinates {
(2,  0.845299462)
(3,  0.845299462)
(4,  0.845299462)
(5,  0.845299462)
(6,  0.845299462)
(7,  0.845299462)
(8,  0.845299462)
(9,  0.845299462)
(10, 0.845299462)
(11, 0.845299462)
(12, 0.845299462)
(13, 0.845299462)
(14, 0.845299462)
(15, 0.845299462)
(16, 0.845299462)
(17, 0.845299462)
(18, 0.845299462)
(19, 0.845299462)
(20, 0.845299462)
(21, 0.845299462)
(22, 0.845299462)
(23, 0.845299462)
(24, 0.845299462)
(25, 0.845299462)
(26, 0.845299462)
(27, 0.845299462)
(28, 0.845299462)
(29, 0.845299462)
(30, 0.845299462)
(31, 0.845299462)
(32, 0.845299462)
(33, 0.845299462)
(34, 0.845299462)
(35, 0.845299462)
(36, 0.845299462)
(37, 0.845299462)
(38, 0.845299462)
(39, 0.845299462)
(40, 0.845299462)
(41, 0.845299462)
(42, 0.845299462)};

        \begin{scope}[shift={(270,90)}]
	\draw[dotted] (0,0) --
		plot[black, dotted] (15,0) -- (30,0)
		node[right]{$2 - 2\sqrt{\frac{1}{3}}$};

        \draw[yshift=1.25\baselineskip] (0,0) --
		plot[mark=*, mark options={fill=white}] (15,0) -- (30,0)
		node[right]{$n B_n^2$};
        \end{scope}

\end{axis}
\end{tikzpicture}
\end{center}

The dotted line on the plot above corresponds to vectors $a, b, c, d$ forming a regular tetrahedron. Let's take $n = 4$ and vectors
\begin{equation*}
\begin{split}
a_1 = \frac{1}{2\sqrt{3}} (1, 1, 1), \; a_2 &= \frac{1}{2\sqrt{3}} (-1, -1, 1), \\
a_3 = \frac{1}{2\sqrt{3}} (-1, 1, -1), \; &a_4 = \frac{1}{2\sqrt{3}} (1, -1, -1).
\end{split}
\end{equation*}

Applying the inverse Hopf map, we obtain
\begin{equation*}
A = \frac{1}{2\sqrt[4]{3}} \; W \left( \begin{matrix}
\sqrt{\sqrt{3} + 1} &\frac{1-i}{\sqrt{\sqrt{3} + 1}} \\
\sqrt{\sqrt{3} + 1} &\frac{-1+i}{\sqrt{\sqrt{3} + 1}} \\
\sqrt{\sqrt{3} - 1} &\frac{1+i}{\sqrt{\sqrt{3} - 1}} \\
\sqrt{\sqrt{3} - 1} &\frac{-1-i}{\sqrt{\sqrt{3} - 1}}
\end{matrix} \right),
\end{equation*}
where $W$ is an arbitrary diagonal unitary matrix. All $2 \times 2$ submatrices of matrix $A$ have the same least singular value equal to $\frac{1}{2} \sqrt{2 - 2\sqrt{\frac{1}{3}}}$, which is less than $\frac{1}{2}$ proven to be the sharp lower bound in the real case (see \cite{Nesterenko2023}).

\section{Remarks}

The provided geometrical interpretation of Problem \ref{p1} can be considered as the application of the results by Hausmann and Knutson identifying the spaces of planar and spatial isoperimetric $n$-gons (up to translations and rotations) with the space of $2$-planes in $\mathds{R}^n$ and $\mathds{C}^n$, respectively. We refer the readers to the paper \cite{HK1997} for precise definitions and formulations.

\section{Acknowledgements}

I would like to thank Igor Makhlin, Stanislav Budzinskiy and the authors of the original paper \cite{GTZ1997} for fruitful discussions.

\bibliographystyle{plain}
\bibliography{lit}

\end{document}